\documentclass[12pt]{amsart}
\usepackage{amssymb}

\def\squareforqed{\hbox{\rlap{$\sqcap$}$\sqcup$}}
\def\qed{\ifmmode\squareforqed\else{\unskip\nobreak\hfil
\penalty50\hskip1em\null\nobreak\hfil\squareforqed
\parfillskip=0pt\finalhyphendemerits=0\endgraf}\fi\medskip}

\newcommand{\PSL}{\mathrm{PSL}}

\newcommand{\Tits}{{}^2F_4(2)'}

\title[Explicit embedding of the Tits group in $E_6$]
{An explicit embedding of the Tits group in the compact real form of $E_6$}
\author{Robert A. Wilson} 
\address{
School of Mathematical Sciences, Queen Mary University of London,
Mile End Road, London E1 4NS, UK}
\email{R.A.Wilson@qmul.ac.uk}
\date{First draft 08/08/2012; this version 21/08/2012}

\begin{document}

\begin{abstract}
We construct a set of 
$27\times 27$ unitary matrices which give an explicit embedding of the Tits group
in the compact real form of the Lie group of type $E_6$. A subset gives an 
embedding of $\mathrm{PSL}_2(25)$ in $F_4$.
\end{abstract}

\maketitle

\section{Introduction}
\label{intro}
It is well-known that the Tits simple group ${}^2F_4(2)'$ embeds in $E_6(\mathbb C)$,
but the published proof of this fact by Cohen and Wales
\cite{CW} does not lend itself very easily to obtaining
an explicit embedding. On the other hand,
the $27$-dimensional representation of $\Tits$, written over
$\mathbb Q(i)$, was constructed by Simon Nickerson in his
PhD thesis \cite{SimonN}, and is available from \cite{webatlas}.
The task then is to determine the hermitian form and the cubic form preserved by
this representation. 
By character theory, both forms are unique up to scalar multiplication, so they must be the ones
which define the compact real form of $E_6$.
One would then want to find a basis with respect to which this is
the standard copy of $E_6$.

\section{Strategy}
The best subgroup to use to exhibit the forms would seem to be 
the maximal subgroup $5^2{:}4A_4$. This acts on the $27$-space as the direct sum of
the $3$-dimensional deleted permutation representation of the quotient $A_4$,
and a monomial representation formed from the $24$ non-trivial representations of $5^2$.
Inside this subgroup is an element of order $12$, whose normalizer is $D_{24}$. Thus one can
generate $\Tits$ with the subgroup $5^2{:}4A_4$ and an involution in the outer
half of this $D_{24}$. When correctly scaled, the $24$ non-trivial eigenspaces of
the $5^2$, together with the $3$ non-trivial eigenspaces of $5^2{:}Q_8$, form an
orthonormal basis. Moreover, the terms of the cubic form are severely restricted by
the action of the $5^2$.

\section{Computing a good basis}
We begin by taking Nickerson's representation, which is written over
$\mathbb Z[\frac12,i]$, and reduce all the matrix entries modulo $41$
(so that we have a field with 4th and 5th roots of unity, but no cube roots).
We take $i=9$, and note that the primitive 5th roots of unity are $-4,16,18,10$.
We find the subgroup $5^2{:}4A_4$ using the words in the standard generators $a,b$
given in \cite{webatlas}, that is $a$ and $$c=b^{(abab^2)^3}.$$
Now we find a subgroup $4A_4$ as the centralizer of the involution $(ac)^6$:
this subgroup may be generated by $ac$ and $$d=(a(ac)^6)^5.$$

We next find two elements 
\begin{eqnarray*}
f_1&=&(a(ac)^6)^2,\cr
f_2&=&(f_1)^{ac}
\end{eqnarray*}
of order $5$, generating the $5^2$, and
look for a non-trivial eigenspace, for example as the common nullspace of $f_1+4$
and $f_2+4$. Similarly we find a non-trivial eigenvector of $5^2{:}Q_8$ in the
simultaneous nullspace of $f_1-1$, $f_2-1$ and $d-1$.
The images of our eigenvectors under the group $4A_4$
give us our $27$ basis vectors
(up to signs).

\section{Computing good generators}
We next look for elements in the centraliser in the Tits group of the involution
$(ac)^6$. Using Bray's algorithm, we quickly find the element 
$$e=((ac)^6(ab^2)^{-1}(ac)^6(ab^2))^4,$$
which lies in the outer half of the involution centralizer of shape
$$2^{1+2+1+2+2}S_3.$$ We want to find an element in this involution centralizer
which inverts $(ac)^4$, so we apply the formula:
$$e'=e(ac)^8e(ac)^4e$$
and find that in fact $e'$ is an involution inverting $ac$. From now on we use
exclusively the new generators
$$f_1,f_2,d,ac,e'$$
for $\Tits$.

\section{Adjusting the basis}
Our next task is to adjust the scalars so that vectors in both parts of the basis
have the same norm. To do this, we write $e'$ with respect to the new basis, and
inspect the $3\times 24$ and $24\times 3$ blocks where the two parts of the
basis interact. We find that in each block the non-zero entries are a particular
scalar multiple of $\pm1,\pm 9$, but they are different multiples in the two blocks. 
So we find a suitable
scalar to multiply one of the orbits of basis vectors by, to make these multiples the same.

At this stage, the top row of the matrix contains $19$ non-zero entries, 
which are 25 (twice), 33 (five times), and $8,10,31$ (four times each),
which are
$-8$ times $2, \pm1, \pm 9$. 
It is possible to adjust some of the basis vectors by factors of $\pm 9$
so that the entries become $-8$ times $2$ (twice), times $1$ (nine times) and
times $-1$ (eight times). Other cosmetic changes, such as re-ordering the
coordinates, can be done according to one's personal taste.

\section{Lifting to complex numbers}
Noticing that $-8\equiv 1/5 \bmod 41$ 
we see that this gives an obvious way to lift to the complex numbers, in such a way
that the top row of the matrix of $e'$
has two entries $2/5$ and 17 entries $\pm1/5$, and therefore has norm $1$.
We take $$z=e^{2\pi i/5}$$ as the lift of $16$, so that $10,37,18$ lift to
$z^2,z^3,z^4$ respectively, and also
$7$ lifts to $$\sigma=(1-\sqrt{5})/2=-z-z^4$$ and $35$ lifts to
$$\tau=(1+\sqrt{5})/2=-z^2-z^3.$$ Then $f_1$ and $f_2$ 
act diagonally 
as follows:
\begin{eqnarray*}
f_1&=&\mathrm{diag}(1,1,1;z,z^2,z^3,z^4;1,1,1,1;\cr
&&\qquad\qquad z^3,z,z^4,z^2;z^3,z,z^4,z^2;\cr
&&\qquad\qquad z,z^2,z^3,z^4;z^2,z^4,z,z^3),\cr
f_2&=&\mathrm{diag}(1,1,1;z^4,z^3,z^2,z;z^4,z^3,z^2,z;\cr
&&\qquad\qquad z^3,z,z^4,z^2;z,z^2,z^3,z^4;\cr
&&\qquad\qquad z^3,z,z^4,z^2;1,1,1,1).
\end{eqnarray*}
The elements $d$ and $ac$ permute the six $4$-spaces which are the fixed spaces of
one of the cyclic subgroups of $5^2$, which we have delimited by semicolons above.
Writing
$$I=\begin{pmatrix}1&.&.&.\cr .&1&.&.\cr .&.&1&.\cr .&.&.&1\end{pmatrix},
J=\begin{pmatrix}.&.&.&1\cr .&.&1&.\cr .&1&.&.\cr 1&.&.&.\end{pmatrix},$$
$$K=\begin{pmatrix}.&1&.&.\cr .&.&.&1\cr 1&.&.&.\cr .&.&1&.\end{pmatrix},
L=\begin{pmatrix}.&.&1&.\cr 1&.&.&.\cr .&.&.&1\cr .&1&.&.\end{pmatrix},$$
we have
$$d=\begin{pmatrix}1&.&.&&&&&&\cr
.&-1&.&&&&&&\cr
.&.&-1&&&&&&\cr
&&&-I&.&.&.&.&.\cr
&&&.&-J&.&.&.&.&.\cr
&&&.&.&.&-iI&.&.\cr
&&&.&.&iI&.&.&.\cr
&&&.&.&.&.&.&iL\cr
&&&.&.&.&.&-iK&.\end{pmatrix}$$ and
$$ac=\begin{pmatrix}
.&1&.&&&&&&\cr
.&.&1&&&&&&\cr
1&.&.&&&&&&\cr
&&&.&.&.&.&K&.\cr
&&&.&.&.&.&.&K\cr
&&&K&.&.&.&.&.\cr
&&&.&K&.&.&.&.\cr
&&&.&.&K&.&.&.\cr
&&&.&.&.&K&.&.\end{pmatrix}.$$
Finally, in order to write down $e'$, we define the matrices
$$
A=\begin{pmatrix}0&\tau&\sigma&0\cr \tau&0&0&\sigma\cr \sigma&0&0&\tau\cr
0&\sigma&\tau&0\end{pmatrix},
B=\begin{pmatrix}-1&\tau&\sigma&-1\cr \tau&-1&-1&\sigma\cr \sigma&-1&-1&\tau\cr
-1&\sigma&\tau&-1\end{pmatrix},$$
$$
C=\begin{pmatrix}1&\sigma&\tau&1\cr \sigma&1&1&\tau\cr \tau&1&1&\sigma\cr 
1&\tau&\sigma&1\end{pmatrix},
D=\begin{pmatrix}1&0&0&1\cr0&1&1&0\cr0&1&1&0\cr1&0&0&1\end{pmatrix},
E=\begin{pmatrix}\tau&1&1&\sigma\cr 1&\sigma&\tau&1\cr 1&\tau&\sigma&1\cr
\sigma&1&1&\tau\end{pmatrix},$$
$$
F=\begin{pmatrix}-1&\sigma&\tau&-1\cr \sigma&-1&-1&\tau\cr \tau&-1&-1&\sigma\cr 
-1&\tau&\sigma&-1\end{pmatrix},
G=\begin{pmatrix}\tau&0&0&\sigma\cr 0&\sigma&\tau&0\cr 0&\tau&\sigma&0\cr
\sigma&0&0&\tau\end{pmatrix},$$
and then we have
$$e'=\frac15\begin{pmatrix}
2&2&1&0&-1&1&1&-1&0\cr
2&1&2&-1&0&0&-1&1&1\cr
1&2&2&1&1&-1&0&0&-1\cr
0&-1&1&A&B&A&D&E&F\cr
-1&0&1&B&C&D&G&F&G\cr
1&0&-1&A&D&E&F&A&B\cr
1&-1&0&D&G&F&G&B&C\cr
-1&1&0&E&F&A&B&A&D\cr
0&1&-1&F&G&B&C&D&G
\end{pmatrix}.$$

\section{The permutation representation on $2304$ points}
The generators we have chosen have the property that the subset
$$\{f_1,f_2,ac,e'\}$$ generates the subgroup $\PSL_2(25)$ of index $2304$ in
$\Tits$. It is easy to see that this subgroup fixes the vector $$(1,1,1;0^{24}),$$
and hence we obtain an orbit of $2304$  images of this vector, on which
$\Tits$ acts. This can be used to provide a (computer-aided) proof that our
group really is the Tits group.

\section{The permutation representation on 1755 points}
A calculation with the Meataxe, or GAP, shows that the $1$-space spanned by
the vector
$$(0,0,0;0,0,0,0;i,1,-1,-i;0^{16})$$
has 1755 images under $\Tits$. The point stabilizer is $2^{1+4+4}{:}5{:}4$, 
namely the cenralizer of the involution $d$, and
the quotient of order $4$ (represented for example by the powers of $(ac)^3$)
multiplies the given vector by powers of $i$.
\section{The cubic form}
Dickson's cubic form \cite{Dickson} for $E_6$ has $45$ terms in the $27$ variables, and if
$uvw$ is one of the terms, then the eigenvalues $\lambda_1,\lambda_2,\lambda_3$
of $f_1$ (or $f_2$) on $u,v,w$ respectively must satisfy
 $$\lambda_1\lambda_2\lambda_3=1.$$
Hence any pair $u,v$ of the last 24  basis vectors lies in at most one triple of the cubic form.
Now a small calculation with the Meataxe finds the (unique) $27$-dimensional quotient
of the symmetric square of this $27$-dimensional representation of $\Tits$,
from which one can easily read off which pairs of basis vectors actually occur.

It turns out that  only terms involving one vector from each of the three blocks of eight actually occur,
as well as terms involving one of the first three basis vectors. 
Labelling the coordinates 
$$(-3,-2,-1;1,2,3,4;\ldots;21,22,23,24)$$
for convenience, the triples which occur are images under the monomial group
$5^2{:}4A_4$ of
$$(-3,-2,-1), (-3,1,4), (1,9,17),(1,10,24).$$
It remains to determine the signs of these terms.
One may be specified arbitrarily, and then the rest can be read off from the
above Meataxe calculation.
They can be taken to be
$+$ on $(-3,-2,-1)$ and $(1,10,24)$, and $-$ on
$(-3,1,4)$ and $(1,9,17)$.

One can now identify the $27$ coordinates 
(with some sign changes) 
with those given in Dickson's original
construction \cite{Dickson}, or other sources, such as \cite{FSG}.
One possible labelling of our $27$ coordinates in the 
notation of \cite{FSG} is as follows:
\begin{eqnarray*}
&(0,0'',0';&
-4,1,8,-5;7,6,3,2;\cr
&&-4',6',3',-5';1',7',2',8';\cr
&&-4'',7'',2'',-5'';6'',1'',8'',3'').\end{eqnarray*}
The signs indicate that the first and fourth coordinate in each block of $8$
must be negated.

\section{Embedding $\mathrm{PSL}_2(25)$ in $F_4$}
Finally, the fixed vector $(1,1,1;0^{24})$ 
of the subgroup $\PSL_2(25)$ has as its support the three coordinates corresponding to
one of the terms of Dickson's cubic form, and therefore can be taken 
as the identity element of the exceptional Jordan algebra. Hence
we obtain as a by-product an explicit embedding of $\PSL_2(25)$ in the compact real
form of the Lie group of type  $F_4$.
\section*{Acknowledgement}
I would like to thank Marek Mitros for asking the question how to find an explicit
embeding of $\Tits$ in the compact real form of $E_6$, and for uncovering some
typographical errors in an ealier draft which, if left uncorrected, 
would have rendered the results incorrect and therefore useless.

\end{document}